\documentclass[dvips,preprint]{imsart}

\usepackage{amsthm,amsmath,amsfonts,amssymb}
\usepackage{enumerate}
\usepackage{graphicx}
\usepackage[latin1]{inputenc}  
\usepackage[T1]{fontenc}
\usepackage{bbm}

\usepackage{verbatim} 
\usepackage{xspace} 
\usepackage[numbers]{natbib}

\RequirePackage[OT1]{fontenc}
\RequirePackage{amsthm,amsmath}

\startlocaldefs

\theoremstyle{plain}
\newtheorem{theorem}{Theorem}

\newtheorem{corollary}{Corollary}

\newtheorem{definition}{Definition}

\newtheorem{lemma}{Lemma}
\newtheorem{proposition}{Proposition}
\newtheorem{remark}{Remark}

\newtheorem*{theorem*}{Theorem}
\newtheorem*{definition*}{Definition}
\newtheorem*{lemma*}{Lemma}
\newtheorem*{theoremfr*}{Théorème}
\newtheorem*{definitionfr*}{Définition}
\newtheorem*{lemme*}{Lemme}
\newtheorem*{proposition*}{Proposition}
\newtheorem*{corollaryfr*}{Corollaire}
\newtheorem*{remark*}{Remark}

\numberwithin{equation}{section}

\newcommand{\espace}{\vspace{.4cm}}

\newcommand{\Real}{\mathbb R}
\newcommand{\N}{\mathbb N}
\newcommand{\Z}{\mathbb Z}

\newcommand{\Pro}{\mathbf{P} }

\newcommand{\Prob}{\mathbb P}
\newcommand{\Probc}{\Prob^c}

\newcommand{\eps}{\varepsilon}
\newcommand{\To}{\longrightarrow}
\newcommand{\Inv}{\frac{1}}

\newcommand{\m}{\mathbf{m}} 
\renewcommand{\m}{\nu}

\newcommand{\un}{\mathbbm 1}

\newcommand{\Tr}{\zeta}

\newcommand{\Vr}{V}
\newcommand{\Vrm}{V^-}

\newcommand{\Tri}{\Tr_\infty}

\newcommand{\Pzp}{\Prob_{0+}^c}
\newcommand{\Pzps}{\Prob_{0+}^{c*}}

\newcommand{\F}{\mathcal F}

\newcommand{\de}{\mathrm{d}}

\newcommand{\ccv}{\exp(-\pi / \sqrt 3)}
\newcommand{\overshoot}{m} 

\newcommand{\ccr}{c_{crit}}
\newcommand{\Q}{\mathbf{Q}}

\newcommand{\Cbb}{\Ccal}
\newcommand{\Ccal}{\mathcal{C}}
\newcommand{\Cc}{\Cbb^{0}}
\newcommand{ \var }{a}
\newcommand{\Np} {\mathcal{N}} 

\newcommand{\D}{\mathcal{D}}

\newcommand{\Qbf}{\Q}

\newcommand{\legal}{\stackrel d =}
\newcommand{\Pf}[1]{P^{\uparrow #1}}
\newcommand{\Pun}{\Pro^1}
\newcommand{\Hp}{h}
\newcommand{\Hb}{\overline \Hp}

\newcommand{\citejacobIHP} {\cite{jacob2010}} 

\newcommand{\SDE} {stochastic differential equation\xspace}

\newcommand{\RLPE} {(RLP)\xspace}

\endlocaldefs

\begin{document}

\begin{frontmatter}

\title{Langevin process reflected at a partially elastic boundary I}
\runtitle{Reflected Langevin processes I}

\begin{aug}
\author{\fnms{Emmanuel} \snm{Jacob}\ead[label=e1]{emmanuel.jacob@normalesup.org}}
\address{Université Paris VI --- UPMC \\
Laboratoire de Probabilités et Modèles Aléatoires\\
4 place Jussieu --- Case courrier 188 \\
75252 Paris cedex 05 \\
\printead{e1}}
\affiliation{Laboratoire de Probabilités et Modèles Aléatoires}
\runauthor{E. Jacob}
\end{aug}

\begin{abstract}
Consider a Langevin process, that is an integrated Brownian motion,
constrained to stay in $[0,\infty)$ by a partially elastic boundary
at 0. If the elasticity coefficient of the boundary is greater than
or equal to $\ccr= \exp(-\sqrt \pi /3)$, bounces will not accumulate
in a finite time when the process starts from the origin with
strictly positive velocity. We will show that there
exists then a unique entrance law from the boundary with zero
velocity, despite the immediate accumulation of bounces. This result
of uniqueness is in sharp contrast with the literature on
deterministic second order reflection. Our approach uses certain
properties of real-valued random walks and a notion of spatial
stationarity which may be of independent interest.
\end{abstract}



\end{frontmatter}

\section{Introduction}

In 1905, Einstein has been the first one to develop the theory of Brownian motion,
providing an explanation of the erratic trajectories of particles observed by Brown eighty years earlier.
He considered time scales no smaller than a ``relaxation time'', so that he could suppose the independence
of the displacements of the particles, and proposed a statistical physics approach,
 working on the probability density of the particle, rather than its paths.
 He obtained that this density should satisfy the heat equation, leading to the now usual Brownian motion model.

Three years later, Langevin proposed his own approach. The particles should simply
satisfy the usual equation of motion, stating that their acceleration, multiplied
by their mass, should be equal to the external forces applied to them. The randomness
is only hidden in these forces, which can be decomposed into a deterministic friction term
and a stochastic term, which we would now call white noise.
This leads to the Langevin equation, which is historically the first example of a stochastic equation.
Its solution yields essentially the same behavior as Brownian motion on large time scales.
But on smaller time scales (smaller than the relaxation time), a Langevin process is fundamentally
different from a Brownian motion, since its paths are $C^1$, and is a more accurate model for real particles.
\espace

Today, a well-known and well studied object is the reflected Brownian motion, used to describe
the trajectories of particles constrained to stay in a domain, and in many other applications.
However, there has been only few studies of the reflected Langevin processes yet, and
this paper proposes such a study. For the sake of simplicity, we consider only
the simplest Langevin process. The space is one dimensional, the only external force is a white noise,\footnote{that
is, we consider no friction force. The relaxation time is then infinite.}
and the particle has mass one.
Then, if $x$ is the initial position of the
particle and $u$ its initial velocity, its path $X$ is simply given by
$$ X_t= x + u t + \int_0^t B_s \de s,$$
where $B$ is the standard Brownian motion driving the motion. We call this process (free) Langevin process,
or integrated Brownian motion. A consequent study can be found in Lachal
\cite{Lachal03}. \espace

Further, suppose this particle is constrained to stay in $[0,\infty)$ by a barrier at 0, in such a
way that when the particle hits the barrier with incoming velocity $v<0$, it will instantly
bounce back with velocity $-c v \ge 0$, where $c\ge 0$ is a parameter called elasticity coefficient
or velocity restitution coefficient. When $c=1$ we say the reflection is perfectly elastic,
when $c=0$ it is said totally inelastic.
The modeling of this barrier naturally involves second order reflection, which can be expressed, for the Langevin process, by
 the following second order \SDE:
$$ \text{\RLPE}\qquad \left\{ \begin{array}{ccl} \label{eqmouvementintro}
X_t &=& x + \displaystyle \int_0^t \dot X_s\de s \\
 \\
\dot X_t &=& u + B_t - (1+ c) \sum_{0< s \le t} \dot X_{s-}
\un_{X_s=0} + N_t,
\end{array}
\right. $$
\begin{eqnarray*}
\text{where}\hspace{-.3cm}&B&\text{is the standard Brownian motion driving the motion, standard}\\
&& \text{meaning that it starts from $B_0=0$ and has variance $t$ at time $t$.} \\
&N&\text{is a continuous nondecreasing process starting from } N_0=0, \\
&& \text{increasing only when the process } (X,\dot X) \text{ is at } (0,0), \\
&&\text{in the sense } \un_{(X_t,\dot X_t) \ne (0,0)} \de N_t\equiv 0. \\
&(x,u)&\text{is the initial or starting condition.}
\end{eqnarray*}
The model and these equations will be further discussed in the
preliminaries. The second order reflection for a particle
submitted to a deterministic force already reveals a formidable
complexity. See the paper of Ballard \cite{Ballard} for
relatively recent results, and also those of Bressan in 1960
\cite{Bressan}, Percivale in 1985 \cite{Percivale}, Schatzman in
1998 \cite{Schatzman} for further reference. In particular, an
analytic force implies the existence of a unique solution,
but this may fail even with a $C^\infty$ force. The main aim of this work is
to show that our stochastic model is nicer, in the sense that there
is always a unique solution to \RLPE, in the weak sense. The
particular case of an inelastic reflection $c=0$ has already been
treated in Bertoin \cite{SDE} (see also \cite{reflecting} and
\cite{jacob2010}) -- though in slightly less general settings, as
the possibility of a nonzero term $N$ was not considered.

Our first observation is that when the starting position is
$(x,u)\ne(0,0)$, then Equations \RLPE have a unique maximal solution
$(X,\dot X)$ killed when hitting $(0,0)$. We shall see in
Preliminaries $(1.2)$ that this hitting time is infinite if and only
if the coefficient $c$ is no less than
the critical value $\ccr=\ccv$. In the sequel we restrict the study
to that case and investigate what happens when the starting condition is
$(0,0)$. It may seem an easy question but once again an analogy with
the deterministic equations enlightens the difficulty of the
problem. \espace

We will prove the existence of a unique law of a solution
to \RLPE, that is, of a unique reflected Langevin process
started from $(0,0)$. Its law is obtained as the weak limit of the law of the
reflected Langevin process starting from 0 with a nonzero velocity
$u>0$, when $u$ goes to 0. We also express directly the law of
the reflected Langevin process started from $(0,0)$, and translated at some
random time.

These results may seem similar to those obtained in the inelastic case
in \cite{SDE}.
However, the behavior of the reflected process is very different
when the elasticity coefficient is nonzero, and so is the whole
study. In a forthcoming paper, we will also investigate the subcritical
case $0<c<\ccr$. In that case too, we will prove the existence of a unique
reflected Langevin process, but once again, the qualitative behavior
of the reflected process being fairly different, we will have to use other specific techniques. \espace

The guiding line in this article is to focus on the velocities of
the process at the bouncing times, and we start with the crucial
observation that the sequence of their logarithms forms a random walk.
We first prove a convergence result for this random walk (Corollary
\ref{coroPro_v-P}). Then we translate it to a convergence result for
the reflected process itself (Lemma \ref{lemma}), through which we can
prove our main results (Theorem 1 and 2).
 \espace

The preliminaries start with an informal discussion about the model,
and an insight into the qualitative behavior of the reflected
process. Then starts the rigorous mathematical study, where we show
in particular the phase transition at the critical value
$\ccr=\ccv$. We end the preliminary section with defining a notion
of spatial stationarity, in an abstract context, and giving an
abstract convergence result using this notion (Lemma
\ref{CONVERGENCE_PN}), which will be proved in the Appendix.
Section~\ref{section_Entering} starts with the statement of our two
theorems, both relying on Lemma \ref{lemma}.
Section~\ref{SectionConvShiftProc} uses renewal theory and
Lemma~\ref{CONVERGENCE_PN} to construct a spatially stationary
process and reduce the proof of Lemma \ref{lemma} to that of
Lemma~\ref{lemmaA}. Section~\ref{section_proof_supercritical}
handles this proof in the supercritical case, thanks to an explicit
construction\footnote{These two constructions in particular may be
of independent interest.} of the spatially stationary random walk.
However this construction does not hold in the critical case, and
Section~\ref{section_proof_critical} completes then the proof,
thanks to a disintegration \addtocounter{footnote}{-1}
formula\footnotemark{} for the spatially stationary random walk.

\section{Preliminaries}

 \subsection{Informal discussion on the model}

\subsubsection*{First order reflected Brownian motion}

We start with a few words about the first order reflected Brownian
motion, that is the common reflected Brownian motion.
In this model the path is driven by a standard Brownian motion $B$.
The reflected path started from $x\ge0$ is requested to be a process
evolving in $[0,\infty)$ solution of the following equation:
 $$   X_t=x+B_t+N_t,$$
where $N$ models the push of the barrier. The process $N$ starts from 0 and is requested to be
nondecreasing, continuous, and increasing only when the
particle is at the barrier, in the sense $\un_{X_t>0} \de N_t=0$.
 Its solution is simply given
pathwise by a so-called Skorohod reflection:
$$N_t=0 \vee \sup_{0\le s \le t} \{-x-B_s\}.$$
We stress that the non-reflected process $x+B$ has continuous but
non-derivable paths, and that, at a given time, the law of the
future of this process just depends on its current position. For
these reasons, first order reflection is the most natural way to
model a barrier.

\subsubsection*{Second order reflected Langevin process}

On the contrary, a Langevin process has a well-defined velocity, and its behavior
in the future depends both of its position and its velocity at
present. Therefore first order reflection does not make much sense.
Second order reflection, as described in the introduction,
is the most natural model to consider, leading to \RLPE.
In this model, the push of the barrier
has a direct effect only on the velocity of the particle. This push necessarily
involves a discontinuous part, each jump modeling a bounce, which is
happening when the particle is at 0 with nonzero incoming velocity.
This yields the sum in the equation, indexed by the bouncing times.
The only restrictive assumption we make about this term is that the
velocity restitution coefficient be a constant parameter.

But the barrier push may also include, in full generality, a continuous component. We
should take into account a continuous process $N$ (possibly degenerate), that increases only when
the particle is at $(0,0)$. In the cases when $c=0$ and the external
force is non-positive, the only second order reflected process
starting from initial condition $(0,0)$ is the process staying at 0.
We then have $N_t= -\int_0^t f(s) \de s$. This example should be
enough to illustrate the importance of the term $N$ (which was not
considered in \cite{SDE,jacob2010}). \espace

The existing works on second order reflection reveal that it is a
much more complex equation than first order reflection. In the
general case there is no uniqueness result and no simple expression of
any solution. Consequently, a pathwise approach has \emph{a priori}
no chance to succeed.

\subsubsection*{Second order reflection and transience hypothesis}

At any instant $t$ such that $(X_t, \dot X_t)\ne (0,0)$, there is
locally no bounce (or only one bounce), and there is local pathwise
existence and uniqueness of a solution to \RLPE. Therefore, for a
starting condition $(x,u)\ne (0,0)$, Equations \RLPE yield a unique
strong maximal solution stopped when hitting $(0,0)$. The whole
difficulty of second order reflection is concentrated near the point
$(0,0)$, where the process $N$ is not necessarily constant, and
where an infinite number of bounces occur on a finite time
interval. \espace

Now, suppose that the following ``transience hypothesis'' holds:
whenever the process is not at $(0,0)$, the maximal solution is
defined for all positive times (without hitting $(0,0)$). Then the
only obstruction to the existence of a unique solution is the
starting condition $(0,0)$. Observe, now, that in our model, a
solution to \RLPE cannot stay locally at 0, as
a Brownian motion is almost surely not monotone on any interval.
Therefore a solution has to take off instantly, and will never be
again in 0 with zero velocity. Obviously, in that case $N$ must
be equal to 0.

Now, if $(X_t, \dot X_t)_{t\ge 0}$ is a solution to \RLPE starting
from $(0,0)$ and $\eps$ is small, then the process $(X_{\eps+t},\dot
X_{\eps+t})_{t \ge 0}$ is a solution to \RLPE starting from the
random position $(X_\eps, \dot X_\eps)$, which is near to $(0,0)$.
This may suggest to study the convergence of the law of the
solution starting from $(x,u)$, when $(x,u)$ goes to
$(0,0)$. And indeed, we will prove that these laws have a unique
limit, which is the law of a solution starting from $(0,0)$.
But while a major part of this discussion was relevant for any
second order reflection, we stress that this particular result
of convergence to the solution starting from $(0,0)$ is specific
to our stochastic model. Some deterministic forces may lead to
unexpected behaviors, which we illustrate with two counterexamples. \espace

Consider
the easiest counterexample to uniqueness, when the force $f$ is
a negative constant $f_0$, and the elasticity coefficient $c$ is larger than one.
We can write explicitly all the solutions starting from $(0,0)$.
They are given by the trajectory constantly staying on the barrier, and by the
trajectories $X^{a,d}$, for $d\ge0$, $1\le a < c$, defined by
$$\begin{array}{ll}
\forall  t \in [0,d],& X^{a,d}(t)=0  \\
\forall n\in \Z, t\in[d+c^{n}a, d+c^{n+1}a], &  X^{a,d}(t)= \frac {f_0}2(t-d-c^n a) (c^{n+1}a+d-t).
\end{array}$$
It is easily seen that the solutions $X^{a,0}$ can all be approached by
the solution starting from an initial condition close to,
but different from, $(0,0)$, while this is not the case for the other solutions,
which stay in $(0,0)$ a certain amount of time. Therefore not only the solution
starting from $(x,u)\ne(0,0)$ does not converge to a unique limit trajectory
when $(x,u)$ goes to $(0,0)$, but also the limit trajectories do not
yield all the solutions starting from $(0,0)$.
However, we should say that this counterexample is called ``pathological'' by physicists,
since the elasticity coefficient is larger than one.

A physically more realistic counterexample is given
 by Ballard in \cite{Ballard}, Section 5.3, for $c=0$.
 A close look to it reveals that there are actually not only
 two solutions starting from $(0,0)$, as indicated by the author,
 but an infinite number of them, each one leaving $(0,0)$ instantly.
 In addition, any of these solutions can be approached by the
 solution starting from an initial condition close to $(0,0)$.

 \subsubsection*{The perfectly elastic reflected Langevin process}

Finally, let us observe that the special case $c=1$ is straightforward for our model:
whatever the initial condition, the reflected Langevin process has the same law as the absolute value of a
non-reflected Langevin process. In addition we can use previous works to understand better the reflected process.

Suppose the starting position be 0 and starting velocity be nonzero.
Introduce $\Tr_0=0$ and $\Tr_{n+1}:= \inf \{t>\Tr_n: X_t=0 \}$ for the sequence of
the successive bouncing times, and $\Vr_n= \dot X_{\Tr_n}$
for the sequence of the velocities of the process at these bouncing times.
The results of McKean~\cite{McKean} show that the sequence $(\Tr_n, \Vr_n)_{n\ge0}$
is a homogeneous Markov chain with explicit transition probabilities. Lachal furthers
this study in~\cite{Lachal97} by giving explicit formulas for the law of $(\Tr_n, \Vr_n)$ for a fixed $n$.

Now, suppose on the contrary that the initial condition is $(0,0)$. Then
the reflected process has an infinite number of bounces just after the initial time.
The works of McKean and Lachal still describe the bouncing times and velocities at
these instants, now thanks to two sequences. The first one corresponding
to the successive bounds happening after time 1, the second one to the
successive bounds happening before time 1, counted backwardly.

In a fairly similar manner, in the general case $c> 0$, we are led to consider
a single sequence but indexed by $\Z$. To do this, the first bounce
for which the velocity is greater than 1 will be chosen as the reference, in the sense that
this bounce will have index 0.

\begin{remark}
Wong also studies in~\cite{Wong66, Wong70} the passage times to zero for a certain stationary process, which is obtained from the Langevin process by an exponential change of scale in both time and space. The passage times to zero of this stationary process are closely related to a certain stationary random walk that we will introduce later on. However, this process shall not be confused with the ``stationary Langevin process'' introduced in~\citejacobIHP. The two processes do not seem to be directly related.
\end{remark}

\subsection{The model, preliminary study}

\subsubsection*{Notations}

We use the notation $\Real_+$ for the set of nonnegative real numbers
$[0,\infty)$, and $\Real_+^*$ for the set of positive real numbers $(0,\infty)$.
Introduce $D= (\{0\}\times \Real_+^* )\cup (\Real_+^*\times \Real)$
and $D^0:=D \cup\{(0,0)\}$. Our working space is $\Cbb$,
the space of càdlàg trajectories $(x,\dot x): [0,\infty) \to D^0$,
 which satisfy
$$x(t) = x(0)+ \displaystyle \int_0^t \dot x(s)\de s.$$
This space is endowed with the $\sigma-$algebra generated
 by the coordinate maps and with the topology induced by
 the following injection:
$$ \begin{array}{ccc}
\Cbb &\to& \Real_+ \times \D(\Real_+) \\
(x,\dot x) & \mapsto & \big(x(0), \dot x\big),
\end{array}
$$
where $\D(\Real_+)$ is the space of càdlàg trajectories on $\Real_+$, equipped with Skorohod topology.
We denote by $(X,\dot X)$ the canonical process and by $(\F_t,t\ge0)$ its natural filtration, satisfying the usual conditions of right continuity and
completeness. Besides, by a slight abuse of notation, when we define a probability measure $P$, we also write $P$ for the expectation under this probability measure. When $f$ is a measurable functional and $A$ an event, we also write $P(f,A)$ for the quantity $P(f \un_A)$.
\espace

For any $(x,u) \in D^0$, the second order reflection of the Langevin
process with starting position $x$ and starting velocity $u$ leads
to Equations \RLPE, which we recall here:
$$ (RLP)\qquad \left\{ \begin{array}{ccl} \label{eqmouvement}
X_t &=& x + \displaystyle \int_0^t \dot X_s\de s \\
 \\
\dot X_t &=& u + B_t - (1+c) \sum_{0< s \le t} \dot X_{s-} \un_{X_s=0} + N_t,
\end{array}
\right. $$ where $B$ is the standard Brownian motion driving the motion and
$N$ is requested to be a continuous nondecreasing process starting from $N_0=0$ and increasing
only when $(X_t,\dot X_t)=(0,0)$. A solution is the quadruplet
$(X,\dot X, N, B)$. For any $(x,u)\in D$, there is a unique solution
to Equations $\RLPE$ \emph{killed at the first hitting time of
$(0,0)$ for the process $(X,\dot X)$.} We write $\Prob_{x,u}^c$ for the law
of $(X,\dot X)$, which is a strong Markov process, and whose first
coordinate will be called the killed reflected Langevin process. We
will almost exclusively consider the case when the starting position
is $0$, and write $\Prob_u^c$ for $\Probc_{0,u}$ (with $u>0$).

\espace

Write $\Tr_0=0$ and
$\Tr_{n+1}:= \inf \{t>\Tr_n: X_t=0 \}$ for the sequence of
successive hitting times of zero, and call an \emph{arch} a part of the
path included between two consecutive hitting times of zero. Figure~\ref{fig_arches} below shows
two complete arches and the beginning of a third one. Write also $\Vrm_n$, and $\Vr_n$ for
the speed of the process just before this $n$-th bounce, and for the
speed of the process just after this $n$-th bounce, respectively, so
that we have $\Vr_n= \dot X_{\Tr_n}=- c \dot X_{\Tr_n^-} = - c
\Vrm_n$.
\begin{figure}[!ht]
\centering 
\includegraphics[width=.85\textwidth]{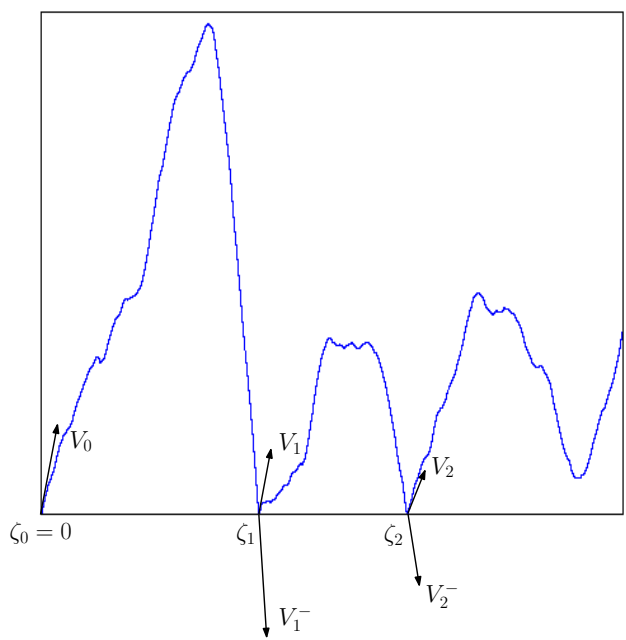}
\caption{First arches of a killed reflected Langevin process}
 \label{fig_arches}
\end{figure}
Please note that the event that for some $n$, we have $\Vr_n=0$, has probability 0. We call time of accumulation of bounces the time $\Tri := \sup(\Tr_n) \in\ (0,\infty]$. It coincides almost surely with the hitting time of $(0,0)$. Next, we will study the sequence $(\Tr_n, \Vr_n)_{n\ge 0}$ and see whether time $\Tri$ is infinite, as in the perfectly elastic case.

\subsubsection*{A phase transition} \label{section_background}

\begin{lemma}\label{lemmePreliminaire}
\begin{enumerate}
\item The law of $( {\Tr_{1}}, {\Vr_{1}} /c)$ under $\Probc_1$ is given by
\begin{equation}\label{loijointe_VrTr}
\begin{array}{l}
 \Inv {\de s \de u} \Prob_1^c\left( ( {\Tr_{1}}, {\Vr_{1}} /c) \in (\de s, \de u)  \right) \\
\hspace{3cm} = \frac {3 u} {\pi \sqrt2 s^2}
\exp(-2\frac{u^2-u+1}{s}) \int_0^{\frac {4u} s} e^{-\frac {3 \theta}
2} \frac {\de \theta} {\sqrt {\pi \theta}}.
\end{array}
\end{equation}

\item Under $\Prob_u^c$, the sequence $\left( \dfrac {\Tr_{n+1}- \Tr_n} {\Vr_n^2}, \dfrac {\Vr_{n+1}} {\Vr_n} \right)_{n\ge0}$ is i.i.d. The common law of its marginals, also independent of $u$, is that of $( {\Tr_{1}}, {\Vr_{1}})$ under $\Probc_1$.

\item In particular, the sequence $\ln(\Vr_n)$ is a random walk. The density of its step distribution $\ln(\Vr_1 /\Vr_0)$ under $\Prob_u^c$ does not depend on $u$ and is given by:
\begin{equation} \label{loi_Vrr}
\Inv {\de u} \Prob_1^c(\ln(\Vr_1) \in \de u)= \frac 3 {2\pi} \frac
{e^{\frac 5 2(u-\ln c)}} {1+e^{3(u-\ln c)}} \de u.
\end{equation}
In particular $\ln(\Vr_1)$ has finite variance and expectation
$$\Prob_1^c(\ln \Vr_1) = \frac \pi {\sqrt 3} +
\ln c.$$

\item
We have, when $t \to \infty$,  
\begin{equation} \label{queue_Trr} \Prob_1^c(\Tr_1>t) \sim c' t^{-\Inv
4},
\end{equation}
where $c'=3 \Gamma(1/4) / (2^{3/4} \pi^{3/2}).$
\end{enumerate}
\end{lemma}

\begin{proof} 
The three first points are essentially results given by McKean~\cite{McKean} or direct consequences of these. The last point is similar to a result of Goldman for the law of the process with zero starting velocity and nonzero starting position~\cite{Goldman2}, and follows from~\eqref{loijointe_VrTr} by standard integral calculus.

For the convenience of the reader, we explain the second point. It follows from the observation that the
variable $(\Tr_n - \Tr_{n-1})/{(\Vr_{n-1})^2}$ (resp. $\Vr_n /
\Vr_{n-1}$) is equal to the duration of the $n-$th arch renormalized
to start with speed one (resp. to the absolute value of the speed of
the process just before its return time to zero, for this
renormalized arch). More precisely:

Recall that, conditionally on $\Vr_n=u$, the process
$(X_{(t+\Tr_n)\wedge\Tr_{n+1}})_{t\ge 0}$ is independent of
$(X_{t\wedge\Tr_n})_{t\ge 0}$ and has the same law as
$(X_{t\wedge\Tr_1})_{t\ge 0}$ under $\Prob_u^c$, thus $(\Tr_{n+1} -
\Tr_n,\Vr_{n+1} /c)$ is independent of $(\Tr_k,\Vr_k)_{k \le n}$ has
the same law as $(\Tr_1,\Inv c \Vr_1)$ under
$\Prob_u^c$. 
 It follows that the variable $(({\Tr_{n+1} -
\Tr_n})/{(\Vr_n)^2},\Vr_{n+1}/ \Vr_n)$ is independent of
$(\Tr_k,\Vr_k)_{k \le n}$ and has the same law as $(\Tr_1, \Vr_1)$
under $\Prob_{1}^c$ (conditionally on $\Vr_n=u$, but this
conditioning can simply be removed). The statement follows.
\end{proof}
From this Lemma we deduce the phase transition phenomena:\\

\begin{corollary} \label{coroTri}
 The time of accumulation of bounces $\Tri$ is:\vspace{.2cm}

 \qquad \begin{minipage}{.8 \textwidth}
    finite $\Prob_u^c-$almost surely if $c < \ccv$, \\
    infinite  $\Prob_u^c-$almost surely if $c \ge \ccv$.
 \end{minipage}\vspace{.2cm}

We thus call $\ccr:=\ccv$ the critical elasticity coefficient. We
call the case $c>\ccr$ the supercritical regime, the case $c<\ccr$
the subcritical regime, the case $c=\ccr$ the critical regime.
\end{corollary}

\begin{proof}
We may express $\Tri$ as the series:
\begin{eqnarray*}
\Tri &=& \sum_{n=1}^\infty \dfrac {\Tr_n - \Tr_{n-1}}
{(\Vr_{n-1})^2} (\Vr_{n-1})^2.
\end{eqnarray*}
For $c<\ccv$, the law of large numbers tells that the sequence
$\Inv k \ln(\Vr_k)$ converges to $\ln(c) + \pi/\sqrt 3<0$ a.s.\ On
the other hand, it follows from \eqref{queue_Trr} that the
expectation of $(\ln(\Tr_{1}))^2$ is finite\footnote{This result was
also stressed by McKean in \cite{McKean}}. Thus, for any fixed $\eps
>0$ there are a.s.\ only a finite number of $k$ such that $\ln((\Tr_k -
\Tr_{k-1})/{(\Vr_{k-1})^2})$ is larger than $\eps k$. We deduce an
a.\thinspace s. exponential decay for the variables $\Tr_{k+1}-\Tr_k$. \emph{A
fortiori} $\Tri$ is a.\,s. finite.

Take now $c \ge \ccv$. For $c>\ccv$, the random walk $\ln \Vr_n$ has a positive drift and is transient.
Thus the sequence $\Vr_n$ is diverging to $+\infty$. As $(\Tr_n -
\Tr_{n-1})/{(\Vr_{n-1})^2}$ is independent of $\Vr_{n-1}$ and has a fixed distribution, we deduce
that $\Tri$ is infinite. For $c= \ccv$, the step distribution has
zero expectation and finite variance, thus the random walk is
recurrent (from the central limit theorem). Then the sequence $\Vr_n$ is recurrent, but it is still not
converging to zero, which is enough to conclude in the same way that $\Tri$ is infinite.
\end{proof}

From now, we will restrict the study to the supercritical and critical regimes, or $c\ge \ccr$, when the transience hypothesis holds.

\subsection{Spatial stationarity} \label{sectionstationarite}

After these first results on the Langevin process, we give the
abstract context for a notion of spatial stationarity and an
important lemma that we will need later.

Write $\Omega$ for the set of sequences indexed by $\Z$
with values in $[-\infty,\infty) \times \Cc$, where $\Cc$ is a
topological space with an isolated point $\emptyset$. For now,
just consider this space as playing an accessory role that will be
clarified later. The set $\Omega$ is endowed with the usual product topology.
An element of $\Omega$ will be written alternatively $\omega$, $(\omega_n)_{n\in \Z}$ or $(\omega_n^1,\omega_n^2)_{n\in \Z}$.

For any real number $x$ we write $T_x$ for the hitting time of
$(x,\infty)$ by the first coordinate, that is
$$ T_x=T_x(\omega) = \inf \{n\in\Z, \omega_n^1>x\}.$$
Under all the measures $\mathrm{P}$ that we will consider on
$\Omega$ we will have
$$ \lim_{-\infty} \omega_n^1 = - \infty, \ \limsup_
{+\infty} \omega_n^1 = + \infty \qquad \mathrm P \text{-almost
surely},$$ and as a consequence $T_x$ will have values in $\Z,$
$\mathrm P$-almost surely. Write $\Omega_0$ for the subset of $\Omega$ consisting in sequences for which $T_0=0$.
Define a spatial translation operator $\Theta$ on $\Omega$, by:
\begin{equation}\label{defTheta}
    \Theta_x(\omega) := (\omega_{n+T_x}^1 - x, \omega_{n+T_x}^2)_{n \in
    \Z}.
\end{equation}
Observe that the range of $\Theta_x$ is always $\Omega_0$, and that
the restriction of $\Theta_0$ to $\Omega_0$ is the identity.
This definition immediately yields a notion of spatial stationarity
for probability laws on $\Omega$:
\begin{definition}\label{defstatiospatiale}
We say that a probability $P$ on $\Omega$ is spatially stationary if
$P \circ \Theta_x=P$ for any $x \in \Real$.
\end{definition}
We also write
$$ \Omega_+ := \left\{ \omega \in \Omega: (\omega_n^1,\omega_n^2)=(-\infty,\emptyset) \text{ for all }n<0 \right\}.$$
An element of $\Omega_+$ shall be thought of as a sequence indexed by $\N$. We write
$\omega^+ \in \Omega_+$ for the projection of $\omega\in \Omega$
defined by:
$$ \omega_n^+ = \left\{\begin{array}{cl}
(-\infty,\emptyset) \qquad &\text{if }n<0, \\
 \omega_n \qquad &\text{if } n\ge0.
\end{array} \right.$$
If $P$ is a probability law on $\Omega$, we write $P_+$ for the image
probability law on $\Omega_+$ by this projection.
Finally the notation $\rightarrow$ simply denotes the weak convergence for probability laws on the topological space $\Omega$.
 The following lemma formulates how convergence results on $\Omega_+$
can imply a convergence result to a
spatially stationary probability measure on $\Omega$.
\begin{lemma} \label{CONVERGENCE_PN}
Let $(P_v)_{v>0}$ be a family of probability laws on $\Omega$. We
suppose that there is a probability law $Q$ on $\Omega_+$ such that:
 \begin{equation} \label{hypothese_CONVERGENCE_PN}
 \forall x \in \Real, \quad (P_v \circ \Theta_x)_+ \rightarrow_{v\to 0} Q.
 \end{equation}
Then there exists a unique spatially stationary probability law $P$ on
$\Omega$ such that $P_+=Q$. Moreover, we have
$$ P_v \circ \Theta_x \rightarrow P.$$
\end{lemma}

The proof of this technical lemma is based on the Kolmogorov
existence theorem. We postpone it to the appendix.

\section{Entering with zero velocity} \label{section_Entering}

Recall that we are in the critical or supercritical regime, $c\ge \ccr$.
%
Write $(S_n)_{n\ge0}$ for the sequence of the logarithms of the
(outgoing) velocity at the successive bounces, defined by
$S_n=\ln(\Vr_n)$. From Lemma~\ref{lemmePreliminaire}, under
$\Probc_{u}$, it is a random walk with step distribution given by
\eqref{loi_Vrr} and drift
$$\mu:= \Prob_1^c(S_1-S_0)= \frac \pi {\sqrt 3} + \ln c = \ln (c/\ccr).$$
In the supercritical case $c>\ccr$ the drift is strictly positive, while  in the critical case $c= \ccr$ the step distribution has zero drift and finite variance.

We introduce the (strictly) ascending ladder height process $(H_n)_{n\ge0}$ associated to the random walk $(S_n)_{n\ge 0}$, that is the random walk with positive jumps defined by $H_0=S_0$ and $H_k = S_{n_k}$, where $n_0=0$ and $n_k=\inf\{n>n_{k-1},S_n>S_{n_{k-1}}\} \in \N$.
In both cases (positive drift, or null drift and finite variance), it is known (see Theorem~3.4 in Spitzer \cite{Spitzer})
that the expectation of the step distribution of $(H_n)_{n\ge0}$, that is $\mu_H:= \Prob_1^c (H_1-H_0),$ belongs to $(0,\infty).$
The probability law
\begin{equation}\label{def_overshoot}
\overshoot(\de y):= \Inv {\mu_H} \Prob_1^c(H_1-H_0 >y) \de y.
\end{equation}
 is known in renewal theory as the stationary law of the overshoot (see also Part~\ref{SectionConvShiftProc}).

We now state our main theorems. The first one is a convergence
result for the probability laws $(\Prob_{u}^c)_{u>0}$ when $u \to
0^+$, while the second one states the weak existence and uniqueness
of solutions to Equations \RLPE with initial condition $X_0=\dot
X_0=0$.
\begin{theorem} \label{theorem}
The family of probability measures $(\Prob_{u}^c)_{u>0}$ on $\Cbb$
has a weak limit when $u \to 0^+$, which we denote by
$\Prob_{0^+}^c$. More precisely, write $\tau_v$ for the instant of
the first bounce with speed greater than $v$, that is
$\tau_v:=\inf\{t>0, X_t=0, \dot X_t > v\}.$ Then the law
$\Prob_{0^+}^c$ satisfies the following conditions:

$$\begin{array}{ll}
(*) &  \left\{\begin{array}{l}
\displaystyle \lim_{v \to 0^+ } \tau_v = 0. \\
\text{For any }u, v>0 \text{, and conditionally on }\dot X_{\tau_v} = u\text{, the process }\\
(X_{\tau_v+ t}, \dot X_{\tau_v+ t})_{t\ge0} \text{ is independent of
}
 (X_s, \dot X_s)_{s< \tau_v} \text{ and has law }\Prob_u^c.
\end{array}
 \right. \\
 \\
(**) & \text{For any }v>0, \text{ the law of }\ln(\dot X_{\tau_v}
/v) \text{ is } \overshoot.
 \end{array} $$
%
\end{theorem}

\begin{theorem} \label{corollary}
\noindent $\bullet$ Consider $(X, \dot X)$ a process of law $\Prob_{0^+}^c$. Then the jumps of $\dot X$ on any finite interval are summable and the process $B$ defined by
$$B_t=\dot X_t + (1+c) \sum_{0<s\le t} \dot X_{s-} \un_{X_s=0}$$
is a Brownian motion. As a consequence the quadruplet $(X,\dot
X,0,B)$ is a solution to $\RLPE$ with initial condition $(0,0)$.

\noindent $\bullet$ For any solution $(X,\dot X, N, B)$ to $\RLPE$
with initial condition $(0,0)$, the law of $(X, \dot X)$ is
$\Prob_{0^+}^c$, and $N\equiv 0$ almost surely.
\end{theorem}

Let us introduce a slightly larger working space, $$\Cbb^*:=\{(x_t, \dot x_t)_{t>0}, \forall \eps>0,(x_{\eps+t}, \dot x_{\eps+t})_{t\ge0} \in \Cbb \}.$$
 We mention that $\Cbb$ can be seen as a subspace of $\Cbb^*$, by removing time 0 from the trajectories. This inclusion is strict: an element of $\Cbb^*$ is a trajectory (indexed by $\Real_+^*$) which does not necessarily have a limit at $0+$.
Both theorems will actually follow from the following lemma, which can be seen as a weak version of Theorem~\ref{theorem}, and whose proof is reported to later.

\begin{lemma} \label{lemma}
There exists a law $\Pzps$ on $\Cbb^*$ such that:

\noindent $\bullet$ We have $\tau_v>0$ for any $v>0$, $\Pzps$-almost
surely.

\noindent $\bullet$ conditions $(*)$ and $(**)$ are satisfied

\noindent $\bullet$  For any $v>0$, the joint law of $\tau_v$ and
$(X_{\tau_v+t}, \dot X_{\tau_v + t})_{t\ge 0}$ under $\Probc_u$
converges weakly, when $u$ goes to 0, to that under $\Pzps$.
\end{lemma}
\begin{proof}[Proof of Theorem~\ref{theorem} and Corollary~\ref{corollary}]
Consider,  under $\Pzps$, the canonical process $(X_t, \dot X_t)_{t>0}$.
From conditions  $(*)$ and the Markov property, we deduce that $(X_t, \dot X_t)_{t>0}$ is a strong Markov process with values in $D$ and transitions that of the reflected Langevin process.

It follows that for any $r>0$, there exists a Brownian motion $B^r$ independent of $\F_r$ and such that, for $t\ge r$,
$$\left\{ \begin{array}{ccl}
X_t &=& X_r + \displaystyle \int_r^t \dot X_s\de s \\
\dot X_t &=& \dot X_r + B^r_{t-r} - (1+c) \sum_{r< s \le t} \dot X_{s-} \un_{X_s=0}.
\end{array}
\right.
$$ 
The Brownian motions $B^r$ are linked by $B^{r}_{t-r}= B^q_{t-q}-B^q_{r-q}$ for $q\le r \le t$.
We introduce $M_s=B_s^{1-s}, \ 0\le s<1$. For any $t<1$, we have
$$ (M_s)_{0 \le s \le t} = (B_t^{1-t} - B_{t-s}^{1-t})_{0 \le s \le t}.$$
Therefore $(M_s)_{0 \le s \le t}$ is a Brownian motion. It follows that $(M_s)_{0\le s <1}$ is a Brownian motion. Write $M_1$ for its limit when $s$ tends to 1.
Now, define the process $B$ by
$$ B_s = \left\{ \begin{array}{ll}
M_1 - M_{1-s} \ &,\ 0 \le s < 1.\\
M_1 + B_{s-1}^1 &,\ 1 \le s.
\end{array} \right.
$$
It is easy to check that $B$ is a Brownian motion and satisfies $B_t - B_r = B^r_{t-r}$ for $t\ge r$. Hence, for $t \ge r$,
%
\begin{equation} \label{eqpreuvecorollaire}
\left\{ \begin{array}{ccl}
X_t &=& X_r + \displaystyle \int_r^t \dot X_s\de s \\
\dot X_t &=& \dot X_r +  B_t-B_r - (1+c) \sum_{r< s \le t} \dot X_{s-} \un_{X_s=0}.
\end{array}
\right.
 \end{equation}

The increments of $\dot X$ are equal to the sum of two terms, on the
one side the increments of $B$, and on the other side, the jumps,
which are happening at the bouncing times. Besides, conditions $(*)$
imply $\dot X_t \un_{X_t=0} \underset{t\to0} \to 0$. That is, the
value of $\dot X$ at a bouncing time is going to 0 when this time
goes to 0. It follows $\dot X_t \underset{t\to0} \to 0$. Therefore
we also have $X_t \to 0$. Consequently, by setting $X_0 = \dot
X_0=0$, we define a process in $\Cbb$. We call its law $\Pzp$. Now,
take again System \eqref{eqpreuvecorollaire} and let $r$ go to 0.
First, we obtain that the sum of the jumps happening just after the
initial time (or in a  finite time interval) is finite. Then we
deduce that under $\Pzp$, $(X, \dot X, 0, B)$ is a solution to
$\RLPE$ with starting condition $(0,0)$. \espace

In summary, we defined a law $\Pzp$ on $\Cbb$ satisfying conditions
$(*)$ and $(**)$, and thus $\tau_v>0$ and $\tau_v \to 0$ almost
surely. Besides, the joint law of $\tau_v$ and $(X_{\tau_v+t}, \dot
X_{\tau_v+t})_{t\ge 0}$ under $\Probc_u$ converges weakly to that
under $\Pzp$. In order to deduce the convergence of $\Probc_u$ to
$\Pzp$, we just need to control what happens on $[0,\tau_v[$. More
precisely, it is enough to control the velocity $\dot X$. Let us
call $M_v$ the supremum of $\dot X_t$ on $[0,\tau_v[$. It will be
enough to prove that when $v$ is small, the variable $M_v$ is small
with high probability, uniformly on $u$ small, in the following
sense:
$$ \forall \eps>0, \forall \delta>0, \exists u_0, v_0>0, \forall 0< u\le u_0, \forall 0<v\le v_0,\quad \Probc_u( M_v \ge \delta)\le \eps.$$
Start from the basic observation $M_v \le v+ \sup_{s,t \in
[0,\tau_v[} |B_t - B_s|$, where $B$ is the underlying Brownian
motion. It follows
\begin{eqnarray*}
\Probc_u(M_v \ge v+ \delta) &\le& \Probc_u(\tau_v \ge \eta) + \Probc_u\left(\sup_{s,t \in [0,\eta)} |B_t - B_s| \ge \delta\right) \\
   &\le&  \Probc_u(\tau_v \ge \eta) + \eps,
\end{eqnarray*}
for a well-chosen $\eta>0$, independent of $u$. Now, by writing the
right side in the form $\Probc_u(\tau_v \ge \eta) - \Pzp(\tau_v \ge
\eta) + \Pzp(\tau_v \ge \eta) + \eps$, and using $\tau_v \underset
{v\to 0} \to 0$, $\Pzp-$a.s., we get that the following inequality
\begin{eqnarray*}
\Probc_u(M_{v} \ge u+ \delta)    &\le&   \Probc_u(\tau_{v} \ge \eta)
- \Pzp(\tau_{v} \ge \eta) + 2 \eps
\end{eqnarray*}
is satisfied for $v$ small enough. Choose $v_0$, smaller than
$\delta$, such that the inequality is satisfied. Then, from the
convergence of the law of $\tau_{v_0}$ under $\Probc_u$ to that
under $\Pzp$, we get that for $u$ smaller than some $u_0>0$, we have
$$  \Probc_u(M_{v_0} \ge 2 \delta) \le 3 \eps.$$
Now it is clear that the inequality stays satisfied for $v<v_0$,
which ends the proof.
The law $\Probc_u$ converges weakly to $\Pzp$, and
Theorem~\ref{theorem} is proved.

\espace Finally, we should prove the uniqueness in
Theorem~\ref{corollary}. Consider any solution $(X, \dot X, N, B)$
to $\RLPE$ with starting condition $(0,0)$. As discussed in the
preliminaries, almost surely, $N=0$ and we have $(X_t,\dot X_t)\ne
(0,0)$ for any positive $t$.
 If the first coordinate $X$ were not coming back to zero at small times,
then there wouldn't be any jumps for $\dot X$ at small times, thus
$X$ would behave like a Langevin process. But this is not possible
as the Langevin process starting from zero with zero velocity does
come back at zero at arbitrary small times. As a consequence, the
process
$(X, \dot X)$ necessarily satisfies condition $(*)$. 
Now, the process $(X_{\tau_v+t}, \dot X_{\tau_v+t})_{t\ge0}$
converges in law to $(X, \dot X)$, thus the law of $(X, \dot X)$ is
an accumulation point of the family $(\Prob_u^c)_{u>0}$ when $u \to
0$. It must coincide with $\Prob_{0^+}^c$.
\end{proof}

The rest of the section is devoted to the proof of
Lemma~\ref{lemma}. It can be sketched as follows. First, using
renewal theory, we get, for any fixed $v>0$, the convergence of the
law of the process $(\dot X_{\tau_v+t})_{t\ge0}$ to a law that can
be described in a simple way. Then Lemma~\ref{CONVERGENCE_PN}
allows, in a certain sense, to include negative times in this
convergence result. The last step will be to prove that $\tau_v$
converges in law to a finite valued random variable.


\subsection{Convergence of shifted processes}\label{SectionConvShiftProc}

We recall the notation $\Vr_n$ for the (outgoing) velocity at the
$n$-th bounce and $S_n$ for its logarithm, for $n\ge0$.
We also
write $\Np_n$ for the translated velocity path starting at the
$n$-th bounce and renormalized so as to start with speed one. That
is, $\Np_n$ is defined by
\begin{equation} \label{N_n}
(\Np_n(t))_{t\ge0} := (\Vr_n^{-1} \dot X (\Tr_n+ \Vr_n^2 t))_{t\ge0}.
\end{equation}
The process $\Np_n$ is independent of $(\dot X_t)_{0\le t \le
\Tr_n}$ and has law $\Prob_1^c$.
The knowledge of the process $X$, or $\dot X$, is equivalent to the knowledge of the sequence $(S_n,\Np_n)_{n\ge0}$, or even just $(S_0, \Np_0)$. But it is more convenient to first prove convergence results about (translations of) the sequence $(S_n,\Np_n)_{n\ge0}$, then deduce results about $X$, which we do.

We work with $\Cc := \Cbb \cup \emptyset$ and we define moreover,
for $n<0$, $(S_n, \Np_n):= (-\infty, \emptyset),$ so that the
sequence $(S,\Np):= (S_n,\Np_n)_{n\in\Z}$ lays  in $\Omega_+$, in
the settings of Section~\ref{sectionstationarite}. We call $\Pro_u$
its law on $\Omega_+$ (or $\Omega$), under $\Prob_u^c$. We also use
the other notations of Section~\ref{sectionstationarite}, such as
$T_x(S)= \inf\{n, S_n\ge x\}$, which we will simply write $T_x$, or
the spatial translation operator $\Theta_x$, defined
by~\eqref{defTheta}. We now aim at establishing convergence results
for the probabilities $\Pro_u\circ \Theta_x$.

First, observe that under $\Pro_u$ and for $n\ge0$,
$(S_{n+1},\Np_{n+1})$ is measurable with respect to $(S_n,\Np_n)$,
and thus $(S,\Np)$ is entirely determined by $(S_0,\Np_0)$, which
follows the law $\delta_{\ln u} \otimes \Prob_1^c$. In other words,
there is a deterministic functional $G$ such that
$(S_n,\Np_n)_{n\ge0}=G(S_0,\Np_0),$ and $\Pro_u$ is the law on
$\Omega$ induced by the law  $\delta_{\ln u} \otimes \Prob_1^c$ for
$(S_0,\Np_0)$. Write now $\Qbf$ for the law on $\Omega_+$ induced by
the law $\overshoot \otimes \Prob_1^c$ for $(S,\Np)$, where the
measure $\overshoot$ is the stationary law of the overshoot we
introduced earlier, defined by \eqref{def_overshoot}.

\begin{lemma}
 For any real number $x$, we have
\begin{equation}
\nonumber (\Pro_u \circ \Theta_x)_+ \Rightarrow_{u \to 0^+} \Qbf
\end{equation}
\end{lemma}

\begin{proof}
Consider the ascending ladder height process $H$ defined at the beginning of Section~\ref{section_Entering}. It is a
random walk with positive jumps and finite expectation. It is nonarithmetic in the sense that its jumping law is not included in $d\Z$ for any $d>0$ (nonarithmeticity is trivial for laws with densities). 
Renewal theory for random walks with positive jumps (see for example
\cite{Gut_SRW}, p.62, or \cite{Feller_Introproba2}, p.355) gives the
following result: the law of the overshoot over a level $x$, that is
$H_{T_x(H)}-x$, converges to $\overshoot$ when $x-H_0$ goes to
infinity. This result is transmitted directly to the random walk
$(S_n)_{n\ge0}$, simply because it has the same overshoot:
$S_{T_x}-x=H_{T_x(H)}-x$.
 Under $\Pro_u$, we have $x-H_0=x- \ln u \to_{u \to 0+} +\infty$.
Hence, when $u$ goes to $0+$, the law of the variable $S_{T_x}-x$
under $\Pro_u$, or, equivalently, that of $S_0$ under $\Pro_u \circ
\Theta_x$, converges to $\overshoot$.

Now, the usual Markov and scaling invariance properties show that
for any $x$, $u$, under $\Pro_u \circ \Theta_x$, $(S_n -
S_0,\Np_n)_{n\ge 0}$ is independent of $S_0$ and has the same law as
$(S_n, \Np_n)_{n\ge0}$ under $\Pro_1$. This altogether establishes
the convergence of $(\Pro_u \circ \Theta_x)_+$ to $\Qbf$.
\end{proof}
Applying Lemma \ref{CONVERGENCE_PN}, we immediately deduce:

\begin{corollary} \label{coroPro_v-P}
For any real number $x$, we have
\begin{equation} \label{Pro_v-P}
 \Pro_u \circ \Theta_x \Rightarrow_{u \to 0^+} \Pro,
\end{equation}
where $\Pro$ is the unique spatially stationary probability measure on
$\Omega$ such that $\Pro_+=\Qbf.$
\end{corollary}

\begin{remark} \label{remark_MASS}
Call $\Pro^1$, resp. $\Qbf^1$, the projection of $\Pro$, resp. $\Qbf$, on the first coordinate. Call $\Theta^1$ the spatial translation operator induced on the first coordinate (defined by $\Theta^1_x(\omega^1) := (\omega^1_{n+T_x} - x)_{n\in \Z}$). Then $\Qbf^1$ is the law of the random walk with starting position distributed according to $\overshoot$. Moreover, we have $\Pro^1_+ = \Qbf^1$, and $\Pro^1$ is spatially stationary. Similar arguments show that $\Pro^1$ is the unique spatially stationary measure such that $\Pro^1_+ = \Qbf^1$. We call it the law of the spatially stationary random walk.
\end{remark}

We now want to deduce Lemma~\ref{lemma} from
Corollary~\ref{coroPro_v-P}. To this end, we have to understand how
to reconstruct $\dot X$ from $\Theta_x(S,\Np)$. We start by working
under $\Pro_u$, for some $u>0$. We introduce an important variable,
$\alpha_x := \tau_{e^x}$, the instant of the first bounce with speed
greater than $\exp(x)$ for the process $(X,\dot X)$.

Observe that the definition \eqref{N_n} of $\Np_n$ induces that the
length of the first arch of $\Np_n$, that is $\Tr_1(\Np_n)$, is
equal to $\Vr_n^{-2}$ times the length of the $(1+n)$-th arch of
$\dot X$.
We may also express $\alpha_x$ as a functional of $ \Theta_x(S,\Np)$ by setting
\begin{equation} \label{alphax}
\alpha_x=  e^{2x} A(\Theta_x(S,\Np)),
\end{equation}
where $A$ is defined by
\begin{equation} \label{A}
 A(\omega) = \sum_{n<0} e^{2 \omega_n^1}
\Tr_1(\omega_n^2),
\end{equation}
 with the convention $\Tr_1(\emptyset)=0$.
Now, the process $(X_t, \dot X_t)_{t\ge \alpha_x}$ is given as the following functional of $ \Theta_x(S,\Np)$:

$$ \left\{ \begin{array}{ccl} \label{reconstructionX_sousPv}
\dot X_t &=& e^{S_{T_x}} \Np_{T_x}(e^{-2 S_{T_x}}(t-\alpha_x))\\
X_t &=& \int_{\alpha_x}^t \dot X_u \de u
\end{array}
\right. \qquad ,\ t\ge \alpha_x. $$

 \espace

 Now, let us work under $\Pro$. It is natural to keep the definition of $\alpha_x$ given by Formula~\eqref{alphax}. Please note however that the sum defining $\alpha_x$ now contains an infinite number of nonzero terms.
\begin{lemma}\label{lemmaA}
1) $\Pro$-almost surely, the time $\alpha_x$ is finite for any $x>0$, and $\alpha_x$ goes to 0 when $x$ goes to $-\infty$,

2) The law of $(\alpha_x,S_{T_x},\Np_{T_x})$ under $\Pro_u$
converges to that under $\Pro$ when $u\to 0^+$.
\end{lemma}
The proof of Lemma~\ref{lemmaA} is postponed to the next subsections. Taking Lemma~\ref{lemmaA} for granted, we may proceed to the proof of Lemma~\ref{lemma}.
\begin{proof}[Proof of Lemma~\ref{lemma}]
 The first part of Lemma~\ref{lemmaA} enables us to define a process $(X_t,\dot X_t)_{t>0}$ on $\Cbb^*$ by
$$ \left\{ \begin{array}{ccl} \label{reconstructionX_sousP}
\dot X_t &=& e^{S_{T_x}} \Np_{T_x}(e^{-2 S_{T_x}}(t-\alpha_x))\\
X_t &=& \int_{\alpha_x}^t \dot X_s \de s
\end{array}
\right. \qquad , \text{ for any } t, x \text{ such that } t\ge \alpha_x. $$
This construction is coherent. We call $\Prob_{0^+}^{c*}$ its law on $\Cbb^*$.

Under $\Prob_{0^+}^{c*}$, the instant $\tau_v:=\alpha_{\ln(v)}$ is
the instant of the first bounce with speed greater than $v$. It is
positive and converges a.s.\ to 0 when $v$ goes to 0. Besides, the
law of $S_{T_{\ln v}} - \ln v$ is equal to $\overshoot$, because by
spatial stationarity, $\Pro \circ \Theta_{\ln v} = \Pro$. Now, take
$x=\ln v$ and $t \ge \alpha_x=\tau_v$ in the formula above. It
follows that under $\Prob_{0^+}^{c*}$, the law of $\ln(\dot
X_{\tau_v} /v)$ is $\overshoot$, and that conditionally on $\dot
X_{\tau_v} = u$, the process $(X_{\tau_v+ t}, \dot X_{\tau_v+
t})_{t\ge0}$ has law $\Probc_u$. We leave to the reader the
verification that it is also independent of $(X_s, \dot X_s)_{0<s<
\tau_v}$. Hence the law $\Pzps$ satisfies conditions $(*)$ and
$(**)$.

The second part of the lemma proves that for any fixed $v>0$, the
joint law of $\tau_v$ and $(X_{\tau_v+t}, \dot
X_{\tau_v+t})_{t\ge0}$ under $\Prob_u^c$ converges to that under
$\Pzps$, as laws on $\Cbb$.
\end{proof}

Finally, all we have to do is to prove Lemma~\ref{lemmaA}. By scaling, it suffices to show that $\alpha_0$ is finite $\Pro-$a.s.\ to prove the first part. We also can suppose $x=0$ for the second part.    
Finally, note that under $\Pro$, we have almost surely $T_0=0$ and hence $\alpha_0= A(\Theta_0(S,\Np))= A(S,\Np)$.

This proof will be based on a more explicit description of the spatially stationary measures $\Pro$ and $\Pro^1$. We must distinguish between the critical and supercritical cases.

\subsection{Proof of Lemma \ref{lemmaA} in the supercritical case} \label{section_proof_supercritical}

%

Throughout this section we suppose that $c>\ccr$. Therefore the drift $\mu= \Pro_1(S_1-S_0)= \frac \pi {\sqrt 3} + \ln c$
is strictly positive. We propose a construction of $\Pro$ based on the introduction of a temporally stationary measure on $\Omega$. If one just considers the first coordinate, this is a construction of the law of the spatially stationary random walk $\Pro^1$, using the temporally stationary random walk.
\espace

First, let us define this temporally stationary random walk. Introduce $P_0$, law of the random walk $(S_n)_{n\in\Z}$ indexed by $\Z$, where $S_0=0$ and $(S_{n+1}-S_n)_{n \in \Z}$ is i.i.d with common law that of the generic step. Then write $P_x$ for the law of $(x+S_n)_{n \in \Z}$ under $P_0$, and set
$$ P_\lambda=\int_\Real P_x \de x.$$
This $\sigma$-finite measure is (temporally) stationary, in the sense that for any $k\in \Z$, the sequences $(S_n)_{n\in\Z}$ and $(S_{k+n})_{n\in\Z}$ have the same law under $P_\lambda$. This term ``law" has to be understood in a generalized sense, that is in settings
where we allow the laws to be not only probability measures but more generally $\sigma$-finite measures. We call this generalized process  of law $P_\lambda$ the (temporally) stationary random walk.

Now start again the same construction, but with adding the second
coordinate. We first recall that under $\Pro_u$ and for $n\ge0$,
$(S_{n+1},\Np_{n+1})$ is measurable with respect to $(S_n,\Np_n)$;
we have $(S_{n+1},\Np_{n+1})=F(S_n,\Np_n)$, where $F$ is a
deterministic functional. For $n\le 0$, consider $\Pi_x^n$ for the
law of $(S_k,\Np_k)_{k\ge n}$, where $\Np_n \legal \Prob_1^c$, $S_n
= x-\ln(\Vr_{-n}(\Np_n))$ (recall that $\Vr_{-n}(\Np_n)$ denotes the
velocity of the particle after the $(-n)-$th bounce), and the
sequence $(S_k,\Np_k)_{k> n}$ is given by $(S_k,\Np_k)=F^{k-n}(S_n,
\Np_n)$.

It should be clear that the laws $\Pi_x^n$, $n\le0$, are compatible. Kolmogorov's existence theorem entails the existence of $\Pi_x$, the law on $\Omega$ under which $(S_k,\Np_k)_{k \ge n}$ has law $\Pi_x^n$ for any $n\le 0$. Then we just define $\Pi_\lambda$ by
$$\Pi_\lambda := \int \Pi_y \de y.$$
Again, this is a $\sigma$-finite (temporally) stationary measure. Besides, the law of the first coordinate $S$ under $\Pi_\lambda$ is $P_\lambda$.

Now, consider the event $\{T_x=n\}$, for $x \in \Real$ and $n\in
\Z$. It should be clear that its measure under $P_\lambda$ is
independent of $x$ and $n$. The following lemma gives its value and
states a link between $\Pi_\lambda$ and $\Pro$, as well as between
$P_\lambda$ and $\Pro^1$ (recall Remark~\ref{remark_MASS} after
Corollary~\ref{coroPro_v-P} for the introduction of the law of the
spatially stationary random walk, $\Pro^1$).

\begin{lemma}\label{LemmaB}
Suppose $c>\ccr$.

1) We have $P_\lambda(T_0=0)=\Pi_\lambda(T_0=0)= \mu \in (0,\infty).$

2) We have $\Pro^1(\cdot)= P_\lambda(\cdot|T_0=0)$ and $\Pro(\cdot)= \Pi_\lambda(\cdot|T_0=0)$.
\end{lemma}

\begin{proof}
Recall that 
$\mu= \Pro_1(S_1-S_0)= \frac \pi {\sqrt 3} + \ln c$
is strictly positive and finite.
 We still write
$(H_n)_{n\ge0}$ for the (strictly) ascending ladder  height process
of the sequence $(S_n)_{n\ge0}$. Its drift $\mu_H= \Pro_1(H_1-H_0)$ is also strictly positive and finite. A result of Woodroofe
\cite{Woodroofe} and Gut \cite{Gut1983} states that, for any $y>0$,
we have
 \begin{equation} \label{gut}
\Inv {\mu_H} P_0 (H_1>y) = \Inv \mu P_0 \left( \inf_{n\ge 1} S_n
> y\right).
\end{equation}
The calculation below follows:
\begin{eqnarray*}
\Pi_\lambda(T_0=0) &=& P_\lambda(T_0=0) \\
&=& \int_0^\infty \de x P_x\left(\sup_{n \le -1} S_n <0\right) \\
&=& \int_0^\infty \de x P_0\left(\inf_{n\ge 1} S_n >x\right) \\
&=& \mu \int_0^\infty \frac {\de x} {\mu_H} P_0(H_1>x) \\
&=& \mu,
\end{eqnarray*}
where we used a symmetry property in the third line.
As $\mu \in (0,\infty)$ we can condition the infinite measure on the event $\{T_0=0\}$ to get
the probability measure
$$\Pi_\lambda(\cdot|T_0=0):=\Inv \mu \Pi_\lambda(\cdot
\un_{T_0=0}).$$ We leave to the reader the simple verification that
this measure on $\Omega$ is spatially stationary in the sense of Definition~\ref{defstatiospatiale} and is projected on
the measure $\Qbf$ on $\Omega_+$. Thus it must coincide with $\Pro$, by Corollary~\ref{coroPro_v-P}.
\end{proof}

We may now prove the first part of Lemma~\ref{lemmaA}.

\begin{proof}[Proof of Lemma~\ref{lemmaA}.1)]
Recall that we need to prove the $\Pro$-a.s.\  finiteness of the sum $A(S,\Np)$.

We start by proving that it is finite $\Pi_x$-almost surely, for a
fixed $x$. Under $\Pi_x$, the sequence $(\Tr_1(\Np_k))_{k\in\Z}$  is
i.i.d with law that of $\Tr_1$ under $\Prob_1^c$. Using the
Borel-Cantelli lemma and estimate \eqref{queue_Trr}, we get that
there are $\Pi_x$-a.s.\ only a finite number of $k>0$ such that
$\Tr_1(\Np_{-k})$ is bigger than $\exp(\sqrt k)$. On the other hand,
the sequence $(S_{-k})_{k\ge0}$ under $\Pi_x$ is a simple random
walk, with an almost sure linear decay. Hence, the sum $A(S,\Np)$ is
finite $\Pi_x$-a.s. It follows that it is also finite
$\Pi_\lambda$-almost surely (by integration) and $\Pro$-almost
surely (by conditioning on a nontrivial event).
\end{proof}
%
%

For Lemma~\ref{lemmaA}.2), we need to prove the weak convergence of
the law of $(\alpha_0, S_{T_0}, \Np_{T_0})$ under $\Pro_u$ to that
under $\Pro$, when $u \to 0^+$. We start by introducing another
notation,
$$\alpha_{x,y}:=\alpha_y - \alpha_x = \sum_{T_x\le n <T_y} \Vr_n^{-2} \Tr_1(\Np_n) \qquad ,\text{ for }x<y.  $$ 
It is clear that under $\Pro$, as well as under $\Pro_u$, we have
almost surely $\alpha_x \underset{x \to -\infty} \to 0$ and
$\alpha_{x,y} \underset{x \to -\infty} \to \alpha_y.$ We also have a
uniform convergence result: the law of the time $\alpha_x$ under
$\Pro_u$ converges in probability to 0 when $x$ goes to $-\infty$,
\emph{uniformly on $u$}, in the following sense:
 \begin{equation} \label{uniformementvers0}
\forall \eps>0, \forall \eta >0, \exists x_0, \forall x\le x_0,
\forall u>0, \quad \Pro_u(\alpha_x \ge \eps) \le \eta.
\end{equation}

Indeed, for any given $\eps>0$ and $\eta>0$, we may choose $y_0$
such that $\overshoot( [0,y_0]) \ge 1-\eta.$ Now, take $u>0$. If
$u>\exp(x)$, then $\alpha_x=0$, and there is nothing to prove. We
suppose $u\le \exp(x).$ From a scaling property, for any $y\ge0$, we
have
\begin{eqnarray*}
\Pro_u(\alpha_x \ge \eps) &=& \Pro_{u e^y} (\alpha_{x+y}\ge \eps e^{2y}) \\
                    &\le& \Pro_{u e^y} (\alpha_{x+y}\ge \eps).
\end{eqnarray*}
Besides, under $\Pro_{u e^y}$, we have $T_{\ln u}=0$ and thus
$\alpha_{x+y}= \alpha_{\ln u, x+y}.$ Hence, we have
\begin{eqnarray*}
\Pro_u(\alpha_x \ge \eps) &\le& \int_{\Real_+} \overshoot(\de y) \Pro_{u e^y} (\alpha_{\ln u, x+y}\ge \eps) \\
&\le& \eta + \int_{[0,y_0]} \overshoot(\de y) \Pro_{u e^y} (\alpha_{\ln u, x+y_0}\ge \eps) \\
&\le& \eta + \int_{\Real_+} \overshoot(\de y) \Pro_{u e^y} (\alpha_{\ln u, x+y_0}\ge \eps) \\
&\le& \eta + \Pro(\alpha_{\ln u, x+y_0} \ge \eps) \\
&\le&  \eta + \Pro(\alpha_{x+y_0} \ge \eps),
\end{eqnarray*}
where the next to last line is a disintegration formula for $\Pro$ at time $T_{\ln u}$ (recall that the law of $S_{T_{\ln u}}-\ln u$ under $\Pro$ is $\overshoot$). 
Now, for $x$ small enough, and uniformly on $u$, we get
$\Pro_u(\alpha_x \ge \eps) \le 2 \eta$. The uniform convergence
result is proved.

We are ready to tackle the proof of Lemma~\ref{lemmaA}.2).
\begin{proof}[Proof of Lemma~\ref{lemmaA}.2)]
It is enough to prove the convergence of the expectation
$\Pro_u( f(S_{T_0},\Np_{T_0}), \alpha_0 \ge a)$ to $\Pro( f(S_{T_0},\Np_{T_0}), \alpha_0 \ge a)$ for any continuous functional $f:\Real \times \Cbb \to [0,1]$ and any $a>0$.

But Corollary~\ref{coroPro_v-P} induces the convergence of the law
of $(\alpha_{x,0},S_{T_0},\Np_{T_0})$ under $\Pro_u$ to that under
$\Pro$. It follows that $\Pro_u(f(S_{T_0},\Np_{T_0}), \alpha_{x,0}
\ge a)$ goes to $\Pro(f(S_{T_0},\Np_{T_0}), \alpha_{x,0} \ge a)$
when $u$ goes to 0. This term in turn converges to
$\Pro(f(S_{T_0},\Np_{T_0}), \alpha_{0} \ge a)$ when $x$ goes to
$-\infty$. As $\alpha_0 \ge \alpha_{x,0}$ for any $x$, it follows
\begin{equation} \label{liminf}
\liminf_{u\to 0} \Pro_u( f(S_{T_0},\Np_{T_0}), \alpha_0 \ge a) \ge
\Pro( f(S_{T_0},\Np_{T_0}), \alpha_0 \ge a).
\end{equation}

On the other hand, for any $\eta>0$, choose $\eps>0$ such that
$\Pro(\alpha_0 \in [a-\eps, a[) \le \eta$, and then choose $x$,
given by the uniform convergence (\ref{uniformementvers0}), such
that for any $u>0$, $\Pro_u(\alpha_x \ge \eps) \le \eta$. Then,
considering the inequality
\begin{eqnarray*}
& &\Pro_u(f(S_{T_0},\Np_{T_0}), \alpha_0 \ge a)\\
&\le &\Pro_u(f(S_{T_0},\Np_{T_0}), \alpha_{x,0} \ge a-\eps) +
\Pro_u(f(S_{T_0},\Np_{T_0}), \alpha_x \ge \eps)
\end{eqnarray*}
and taking the $\limsup$, we get
\begin{eqnarray*}
 \limsup_{u\to 0} \Pro_u(f(S_{T_0},\Np_{T_0}), \alpha_0 \ge a) &\le& \Pro(f(S_{T_0},\Np_{T_0}), \alpha_{x,0} \ge a-\eps) + \eta \\
   &\le& \Pro(f(S_{T_0},\Np_{T_0}), \alpha_0 \ge a) + 2 \eta.
\end{eqnarray*}
This together with \eqref{liminf} gives the desired result.
\end{proof}

\espace

We finish this subsection with a corollary of Lemma~\ref{LemmaB}.
\begin{corollary} \label{-S_-n}
Under $\Pro^1$, conditionally on $S_0=x\ge0$, the sequence $(-S_{-n})_{n \ge 0}$ has the law of the random walk starting from $-x$ and conditioned to stay positive at times $n\ge1$.
\end{corollary}
\begin{proof}
Under $P_\lambda$ and conditionally on $S_0=x$, the sequence $(-S_{-n})_{n\ge 0}$ has the law of the random walk starting from $-x$. 
The event $\{T_0=0\}$, which is also equal to the event $\{S_0>0, \forall n<0, S_n<0\}$, has a positive and finite probability when $x\ge0$. The expression of $\Pro^1$ given in Lemma~\ref{LemmaB} directly implies the corollary.
\end{proof}

\subsection{Proof of Lemma \ref{lemmaA} in the critical case} \label{section_proof_critical}

In the critical case, we certainly can define $P_\lambda$ and $\Pi_\lambda$ as before, but under these measures the time $T_0$ is almost surely equal to $-\infty$. Lemma~\ref{LemmaB} thus fails, and so does the previous construction of $\Pro^1$ and $\Pro$.

However, an analogue of Corollary~\ref{-S_-n} will stay true and induce another construction of the law of the spatially stationary random walk $\Pro^1$. We will then use it to prove again the $\Pro-$almost sure finiteness of $\alpha_0$, and Lemma~\ref{lemmaA} will follow from the same arguments as before. Throughout this subsection we assume $c=\ccr$.

\subsubsection{The spatially stationary random walk in the critical case.}

In order to formulate the analogue of Corollary~\ref{-S_-n}, we need to define the ``random walk conditioned to stay positive'' for a random walk with null drift, for which the event of staying positive for all positive times has probability 0. This is done in \cite{Bertoin_Doney}. We recall it here briefly.

Write as usual $P_x$ for the law of the random walk starting from position $x$.
If you write $(D_n)_{n\ge 0}$ for the strictly descending ladder height process (defined in the exact similar way as the strictly ascending ladder height process, and also equal to the opposite of the strictly ascending ladder height process of $\widehat S:=-S$)
, the renewal function $\Hp$ is defined by
$$ \Hp(x):=\sum_{k=0}^\infty P_x(D_k\ge 0).$$
In particular $\Hp$ is non-decreasing, right-continuous, and we have $\Hp(0)=1$ and $\Hp(x)=0$ for $x<0$.
The renewal function is invariant for the random walk killed as it enters the negative half-line. It enables us to define the process conditioned on never entering $(-\infty,0)$, thanks to a usual $h-$transform, in the sense of Doob. That is, the law of this process starting from $x>0$, written $\Pf 0_x$, is defined by
\begin{equation}
\Pf 0_x(f(S)) = \Inv {\Hp(x)} P_x\Big(f(S) \Hp(S_n),\inf_{k\le n} S_k \ge0\Big)
\end{equation}
for any $f(S)=f(S_0,..., S_n)$ functional of the first $n$ steps. 
For any $a \in \Real$ and $x>a$, we also write $\Pf a_x$ for the law of the random walk starting from $x>a$ and conditioned on never entering $(-\infty,a)$, defined in the exact same way, by
\begin{equation} 
\Pf a_x(f(S)) = \Inv {\Hp(x-a)} P_x\Big(f(S) \Hp(S_n-a),\inf_{k\le n} S_k \ge a\Big)
\end{equation}
for any $f(S)=f(S_0,..., S_n)$ functional of the first $n$ steps.
 The only other thing we will need to know about $\Hp$ is the following sub-additive inequality, which is a consequence of a Markov property:
\begin{equation}\label{H_concave}
\Hp(x+a)-\Hp(x) \le \Hp(a), \qquad x,a>0.
\end{equation}
Recall that $\mu_H$ is the drift of the strictly ascending ladder height process and write $p(x,y)$ for the transition densities of the random walk.
The following proposition gives a disintegration description of the spatially stationary random walk, which is very similar to that of the spatially stationary Lévy process introduced by Bertoin and Savov in \cite{Bertoin-Savov}. 

\begin{proposition}\label{Pun}
The measure $$\m(\de x \de y):=\Inv {\mu_H}p(0,x+y) \un_{x\ge0, y\ge 0} h(x) \de x \de y$$ is a probability law.
The law of $\Pun$ is determined by: \\
$\bullet$  Under $\Pun$,  $(-S_{-1},S_0)$ has the law $\m$. \\
$\bullet$ Conditionally on $-S_{-1}=x$ and $S_0=y$, the processes $(-S_{-n-1})_{n\ge0}$ and $(S_n)_{n\ge0}$ are independent, the law of $(-S_{-n-1})_{n\ge0}$ is $\Pf0_x$, that of $(S_n)_{n\ge0}$ is $P_y$.
\end{proposition}
The measure $\m$ is nothing else than the stationary joint law of the overshoot and the undershoot.
The proof of this proposition will last until the end of the subsection.
As a preliminary, we introduce a crucial though rather simple lemma.
\begin{lemma}\label{lemmaPf0}
For any $0 \le a \le x$, we have:
\begin{eqnarray} \label{Pf0fa1}
\Pf0_x\Big(\inf_{n\ge0} S_n \ge a\Big) &=& \frac {\Hp(x-a)}{\Hp(x)} \\  \label{Pf0fa2}
\Pf0_x\Big(\cdot | \inf_{n\ge0} S_n \ge a\Big)&=& \Pf a_x(\cdot).
\end{eqnarray}
\end{lemma}

\begin{proof}
By expressing the event $\{\inf_{k\ge0} S_k \ge a\}$ as the limit of the events $ \{ \inf_{0\le k \le n} S_k \ge a \},$ we get
\begin{eqnarray*}
  && \Pf0_x\Big(\inf_{0\le k \le n} S_k \ge a\Big) \\
   &=& \Inv {\Hp(x)} P_x\Big(\Hp(S_n),\inf_{0\le k \le n} S_k \ge a\Big)  \\
   &=&  \Inv {\Hp(x)} P_x\Big(\Hp(S_n- a ),\inf_{0\le k \le n} S_k \ge a\Big) \\
   && + \Inv {\Hp(x)} P_x\Big(\Hp(S_n) - \Hp(S_n-a),\inf_{0\le k \le n} S_k \ge a\Big).
\end{eqnarray*}
The first term of the sum is equal to $\frac{\Hp(x-a)}{\Hp(x)}$ because the function $\Hp(\cdot-a)$ is invariant for the random walk killed
when hitting $(-\infty,a)$. The second term is positive and bounded from above by $\frac {\Hp(a)}{\Hp(x)} P_x(\inf_{0\le k \le n} S_k \ge a),$ which goes to 0 when $n$
goes to $+\infty$. This proves equation \eqref{Pf0fa1}. Then \eqref{Pf0fa2} is straightforward: Indeed, for $f(S)=f(S_0,...,S_n)$ functional of the first $n$ steps, we have:
\begin{eqnarray*}
 && \Pf0_x\Big(f(S)|\inf_{k \ge 0} S_k \ge a\Big) \\
 &=&  \Inv {\displaystyle \Pf0_x\Big(\inf_{k \ge 0 } S_k \ge a\Big)}   {\Pf0_x\left(f(S) \Pf0_{S_n}\Big(\inf_{k \ge 0 } S_k \ge a\Big), \inf_{0\le k \le n} S_k \ge a \right)}  \\
   &=&  \frac {\Hp(x)}{\Hp(x-a)}. \Inv {\Hp(x)} P_x\left(f(S) \Hp(S_n) \frac {\Hp(S_n - a)}{\Hp(S_n)}, \inf_{0\le k \le n} S_k \ge a\right) \\
   &=& \Pf a_x(f(S)).
\end{eqnarray*}
\end{proof}
Now, recall that the invariance property of $\Hp$ yields that, for any $x\ge 0$, we have
$$ \Hp(x) = P_x(\Hp(S_1) \un_{S_1\ge 0}).$$
Define $\Hb$ by $\Hb(x) := P_x(\Hp(S_1), S_1\ge 0)$ for \emph{any real number} $x$. Thus for $x\ge0$, $\Hb$ and $\Hp$ coincide, but for $x<0$ they certainly don't. This enables us to define, for any $x, a \in\Real$, the law $\Pf a_x$ of the random walk starting from $x$ and conditioned on never entering $(-\infty, a)$ \emph{at times $n\ge1$}, by the formula:
\begin{equation}\label{Pf0x2}
\Pf a_x(f(S)) = \Inv {\Hb(x-a)} P_x\Big(f(S) \Hp(S_n-a), \inf_{1\le k \le n} S_k \ge a\Big)
\end{equation}
for any functional $f(S)=f(S_0,...,S_n).$
This definition is of course consistent with our previous notations. The following generalization of Lemma~\ref{lemmaPf0} and its corollary are consequences of straightforward calculations, that we leave to the interested reader
\begin{lemma}
For any $y\le a$, any $x \in \Real$, we have
\begin{eqnarray} \label{Pf0fa_generalise}
 \Pf y_x(\inf_{n\ge1} S_n \ge a) &=& \frac {\Hb(x-a)} {\Hb(x-y)} \\
     \Pf y_x(\cdot | \inf_{n\ge1} S_n \ge a)&=&\Pf a_x(\cdot)
\end{eqnarray}
\end{lemma}
\begin{corollary}
Write $\m_-$ (resp. $\m_+$) for the first (resp. second) marginal of $\m$. These measures on $\Real_+$ are given for $x,y>0$, by
\begin{eqnarray*}
 \m_-(\de x) &=& \Inv {\mu_H} \Hp(x) P_0(S_1 \ge x) \de x. \\
 \m_+(\de y) &=& \Inv {\mu_H} \Hb(-y) \de y.
\end{eqnarray*}
Moreover,
\begin{eqnarray*}
 P_{-\m_-} (S_1 \in \de y | S_1 \ge0 ) &=& \m_+(\de y) \\
 \Pf0_{-\m_+}(\de x) &=& \m_-(\de x),
\end{eqnarray*}
where we have written $P_{-\m_-}(...)$ for $\int P_{-x}(...) \m_-(\de x)$, as well as
$P_{-\m_+}(...)$ for $\int P_{-x}(...) \m_+(\de x).$
\end{corollary}

This corollary should make the introduction of the measure $\m$ in the proposition more transparent. Indeed, it gives us two alternative ways of defining the measure $\Pun$. First, take $S_0$ distributed according to $\m_+$ and, conditionally on $S_0=y$, take $(S_n)_{n\ge0}$ of law $P_y$ and $(-S_{-n})_{n\ge0}$ independent and of law $\Pf0_{-y}$ (in the sense defined just before). Second, take $-S_{-1}$ distributed according to $\m_-$ and, conditionally on $S_{-1}=-x$, take $(S_{n-1})_{n\ge0}$ of law $P_{-x}$ 
conditioned on having a first jump no smaller than $x$, and $(-S_{-n-1})_{n\ge0}$ independent and of law $\Pf0_{x}$.

\begin{proof}[Proof of the proposition]

We need to prove three things, the fact that $\m$ is a probability measure (that is, has mass one), the fact that $\Pun$ is spatially
stationary, and the equality $\Pun_+=\Qbf$.
We start with the spatial stationarity. Fix $a>0$. We should prove that $S=(S_n)_{n\in \Z}$ and $R:= \Theta_a(S)= (S_{T_a+n} - a)_{n\in \Z}$ have the same law under $\Pun$.
\espace

We introduce the notation $L_a$ for the instant of the last passage under level $a$ for the process $S$.
Besides, observe that $T_a$ is also equal to the instant of the last passage under level $a$ for the process $(-R_{-n})_{n\ge0}$. \emph{Suppose that we proved} that $((T_a, -R_{-n})_{0 \le n \le T_a})$  has the same law as the process $(L_a, (S_n)_{0 \le n \le L_a})$ under $\Pf0_{- \m_+}$.
Then, conditionally on $-R_{-T_a} = z$, it is clear that the process $(-R_{-n-T_a})_{n\ge0} = (a-S_{-n})_{n\ge0}$ is independent of $(-R_{-n})_{0 \le n \le T_a}$ and follows the law $\Pf a_z$. Besides, for a process $S$ under $\Pf0_{- \m_+}$, conditionally on $S_{L_a}=z$, the process $(S_{n+L_a})_{n\ge0}$ is independent from $(S_n)_{0 \le n \le L_a}$ and follows the law $\Pf a_z$. This altogether proves that the process $(-R_{-n})_{n\ge0}$  follows the law $\Pf0_{- \m_+}$. Finally, from a Markov property, it is clear that given $R_0=y$, the process $(R_n)_{n\ge0}$ is independent of $(R_n)_{n\le 0}$ and follows the law $P_y$, thus the law of $(R_n)_{n\in \Z}$ is $\Pun$.

Therefore, the only thing we still need to prove is the following duality property\footnote{This property also finds its analogue in \cite{Bertoin-Savov}, in their Theorem~2.}: the variable $(T_a, (-R_{-n})_{0 \le n \le T_a})$ has the same law as the variable
$(L_a, (S_n)_{0 \le n \le L_a})$ for a process $S$ of law $\Pf0_{- \m_+}$.
Fix $n\ge0$ and $f:\Real^{n+1} \to \Real$ a positive continuous functional. We should prove the following equality:
$$\Pun(f((-R_{-k})_{0\le k \le n}) \un_{T_a=n}) = \Pf0_{-\m_+}(f((S_k)_{0\le k \le n}) \un_{L_a = n}).$$
The case $n=0$ is particular and follows from this calculation:
\begin{eqnarray*}
  \Pf0_{-\m_+}(-S_0 \in \de x, L_a=0) &=& \Pf0_{-x}(\inf_{k\ge1} S_k \ge a) \m_+( \de x)\\
   &=&  \Inv {\mu_H} \Hb(-x) \frac {\Hb(-a-x)} {\Hb(-x)} \de x \\
   &=& \m_+(a+\de x) = \Pun(R_0 \in \de x, T_a = 0).
\end{eqnarray*}
In the case $n>0$, we write $\tilde f((S_k)_{0\le k \le n}) := f( (a-S_{n-k})_{0 \le k \le n})$, the usual duality property for random walks stating
$$P_x(f(S) \un_{a- S_n \in \de y }) \de x = P_y ( \tilde f(S) \un_{a- S_n \in \de x }) \de y.$$
We are ready to calculate
\begin{eqnarray*}
\Pun(f((-R_{-k})_{0\le k \le n})\un_{T_a=n})  &=& \Pun(f((-R_{-k})_{0\le k \le n})\un_{T_a=n}) \\
&=& \int\int_{\Real_+ \times [0,a)} \m_2(\de x \otimes \de y),
\end{eqnarray*}
where $\m_2(\de x \otimes \de y)$ is equal to
\begin{eqnarray*}
  \! &&\m_+( \de y ) P_y( \tilde f((S_k)_{0\le k \le n}), S_n-a \in \de x, \forall 0\le i <n, S_i \le a)\\
  \! &=&\!\Inv {\mu_H} \Hb(-y) \de x  P_{-x}(f((S_k)_{0\le k \le n}), a-S_n \in \de y, \forall 0< i \le n, S_i \ge 0) \\
  \! &=&\!\frac {\Hb(-y) \de x} {h(a-y)\mu_H} P_{-x}\Big(f((S_k)_{0\le k \le n}) h(S_n), a-S_n \in \de y, \forall 0< i \le n, S_i \ge 0\Big).
\end{eqnarray*}

Using then \eqref{Pf0x2} and \eqref{Pf0fa_generalise}, it follows
\begin{eqnarray*}
&& \Pun(f((-R_{-k})_{0\le k \le n})\un_{T_a=n}) \\
   &=&  \int\int_{\Real_+ \times [0,a)} \!\m_+(\de x) \Pf0_{-x}\big(f((S_k)_{0\le k \le n}), a-S_n \in \de y\big) \Pf0_{a-y}\Big( \inf_{k\ge 1} S_k \ge a\Big)  \\
   &=&  \int_{\Real_+} \m_+(\de x) \Pf0_{-x}\Big(f((S_k)_{0\le k \le n}), S_n <a, \inf_{k >n} S_k \ge a\Big) \\
   &=& \Pf0_{-\m_+} \big(f((S_k)_{0\le k \le n}) \un_{L_a = n}\big).
\end{eqnarray*}
The measure $\Pun$ is thus spatially stationary.
\espace

Now the two facts that $\m$ has mass one and that $\Pun_+ = \Qbf^1$ both follow from the equality
$$ \Hb(-y) = P_0(H_1 \ge y)$$
for $y\ge 0$ (recall that $H$ is the strictly ascending ladder height process).
Fix some $y \ge 0$. We already know from \eqref{Pf0fa_generalise} that $\Hb(-y) = \Pf0_0(\inf_{n\ge0} S_n \ge y)$, thus we should prove
\begin{equation}\label{dualite_Pf_H}
P_0(H_1 \in \de y) = \Pf0_0(\inf_{n\ge0} S_n \in \de y).
\end{equation}
This will be a consequence from another duality argument.
Write $T_{inf}$ for the instant when $S$ hits its minimum on times $n\ge 1$. Write $\tilde T_1:= \inf\{n>0, S_n > S_0\}$ (so that $S_{\tilde T_1} = H_1$).
Then $(S_k)_{0\le k \le \tilde T_1}$ under $P_0$ and $(S_k)_{0 \le k \le T_{inf}}$ under $\Pf0_0$ are in duality.
Indeed, fix $n>0$ and $f(S)= f((S_k)_{0 \le k \le n})$ a positive continuous functional. Write also
$\tilde f ( (S_k)_{0 \le k \le n}) := f( (S_n - S_{n-k})_{0 \le k \le n})$. Then,
\begin{eqnarray*}
&& P_0^{\uparrow 0}(f(S)\un_{\{T_{inf}=n\}}) \\
&=&
P_0^{\uparrow 0}\Big(f(S),\inf_{1\le k\le n-1}S_k>S_n,\inf_{k\ge n+1}S_k\ge S_n\Big)
\\
&=&
P_0^{\uparrow 0}\left(f(S)P_x^{\uparrow 0}\Big(\inf_{k\ge 1}S_k\ge x\Big)
\Big|_{x=S_n},\inf_{1\le k\le n-1}S_k>S_n\right)
\\
&=&
P_0^{\uparrow 0}\left(\frac{f(S)}{h(S_n)},\inf_{1\le k\le n-1}S_k>S_n\right)
\\
&=&
P_0\Big(f(S),\inf_{1\le k\le n-1}S_k>S_n\ge 0\Big)
\\
&=&
P_0\Big(\tilde{f}(S),\sup_{1\le k\le n-1}S_k<0, S_n\ge 0\Big)
\\
&=&
P_0\Big(\tilde{f}(S),\sup_{1\le k\le n-1}S_k<0, S_n>0\Big)
\\
&=&
P_0\big(\tilde{f}(S)\un_{\{\tilde{T}_1=n\}}\big).
\end{eqnarray*}

This duality property implies in particular \eqref{dualite_Pf_H}.
\end{proof}

\subsubsection{Finiteness of $\alpha_0$ in the critical case.}

The only thing we actually need from the last subsection is the fact that under $\Pro^1$ (or, equivalently, under $\Pro$), the sequence $(-S_{-n})_{n\ge1}$ is a random walk conditioned to stay positive, with some initial law. The paper \cite{HKK} gives very precise results about the behavior of this random walk conditioned to stay positive, and we deduce in particular the following rough bounds that are sufficient for our purposes:

\begin{lemma}
For any $\eps>0$, we have
\begin{equation} \label{vitesse_sqrt}
 n^{-\Inv 2 + \eps} S_{-n} \to -\infty
\end{equation}
when $n \to \infty,$ $\Pro$-a.s.
\end{lemma}

We now work under $\Pro$ and we recall that $\alpha_0$ is then given by
$$ \alpha_0= \sum_{n<0} e^{2 S_n} \Tr_1(\Np_n).$$
We write $L_n:=e^{2 S_n} \Tr_1(\Np_n)$ for the duration of the arch of index $n$. We need to transfer the results about the behavior of $(S_{-n})$ to results about the behavior of $(L_{-n})$. This is made possible by the following lemma:

\begin{lemma}
1) Under $\Pro$ and conditionally on a realization $(S_n)_{n\in\Z}=(s_n)_{n\in\Z}$, the variables $(L_n)_{n\in \Z}$ are mutually independent, and the law of $L_n$ is that of $\Tr_1$ under $\Prob_{\exp(s_n)}^c(\cdot|\Vr_1= \exp(s_{n+1}) )$.

2) If $u,v \le \var$ for some real number $\var$, then
\begin{equation}\label{borneTr1}
\Prob_{u}^c(\Tr_1>t \var^2 |\Vr_1= c v ) \le \frac {16 \sqrt 2}{3\sqrt \pi} t^{-\frac 3 2}.
\end{equation}
\end{lemma}
\begin{proof}
The result of the first part is easy for $(L_n)_{n\ge0}$, and we get the result for $(L_n)_{n\in \Z}$ by spatial stationarity.

For the second part, recall that the law of the couple $(\Tr_1,\Vr_1)$ under $\Prob_u^c$ is known (see Lemma~\ref{lemmePreliminaire}, Formulas \eqref{loijointe_VrTr} and \eqref{loi_Vrr}). We obtain, explicitly:
\begin{eqnarray*}
&& \Inv {\de s} \Prob_u^c\left( {\Tr_{1}} \in \de s | {\Vr_{1}} =c v  \right) \\
&& \hspace{2cm} = \frac {\sqrt 2 (u^3 + v^3)} { s^2 u^{\Inv 2} v^{\Inv 2}} \exp\Big(-2\frac{v^2-uv+u^2}{s}\Big) \int_0^{\frac {4 u v} s} e^{-\frac {3 \theta} 2} \frac {\de \theta} {\sqrt {\pi \theta}}.
\end{eqnarray*}
Provided that we take $u, v \le \var$ we get
\begin{eqnarray*}
  \Inv {\de s} \Prob_u^c\left( {\Tr_{1}} \in \de s | {\Vr_{1}} =c v  \right) &\le&
  \frac {2 \sqrt 2 \var^3} { s^2 u^{\Inv 2} v^{\Inv 2}}  \int_0^{\frac {4 u v} s} \frac {\de \theta} {\sqrt {\pi \theta}} \\
  &\le&  \frac{8 \sqrt 2} { \sqrt \pi} \var^3 s^{-\frac 5 2}.
\end{eqnarray*}
Integrating this inequality between $t \var^2$ and $+\infty$ gives \eqref{borneTr1}.
\end{proof}
The $\Pro$-almost sure finiteness of $\alpha_0$ follows straightforwardly. Write
$$A_n=e^{S_n} \vee \frac {e^{S_{n+1}}} c,$$
and, for $n>0$, write $E_n$ for the event $$L_{-n} \ge n \ A_{-n}^2.$$
The lemma states that the probability of $E_n$ is bounded above by a constant times $n^{-\frac32}$. Hence only a finite number of $E_n$ occur, almost surely. This together with \eqref{vitesse_sqrt} gives that the $(L_{-n})_{n\ge0}$ are summable, almost surely.
This shows the $\Pro$-almost sure finiteness of $A(S,\Np)$ and concludes the proof.

\appendix

\section{Proof of Lemma~\ref{CONVERGENCE_PN}}

The uniqueness stated in the lemma is immediate. Indeed, if $P$ and $P'$ are two
probability laws satisfying the conditions of Lemma~\ref{CONVERGENCE_PN}, then we
have $(P\circ \Theta_x)_+ = P_+ = Q= (P' \circ \Theta_x)_+,$ for any real $x$, leading to $P=P'$.
The existence result is based on Kolmogorov's existence theorem, as follows.

First, note that we have $\Theta_x \circ \Theta_y= \Theta_{x+y}$
for any $x, y$ real numbers.
Consider $(P_u)_{u>0}$ and $Q$ satisfying the hypothesis~\eqref{hypothese_CONVERGENCE_PN}.
Our first observation is that $Q$ is necessarily concentrated on $\Omega_0\cap\Omega_+$,
and enjoys already the following ``positive translation invariance property''. Consider any $x>0$.
The equality $ ((P_u \circ \Theta_0)_+ \circ
\Theta_x)_+ = (P_u \circ \Theta_x)_+$ immediately yields, letting $u$ go to 0
and using for each term the hypothesis~\eqref{hypothese_CONVERGENCE_PN}, the equality
\begin{equation} \label{invariance_Q}
(Q\circ \Theta_x)_+=Q.
\end{equation}

For $x_1< ... < x_n$ real numbers, let $Q^{x_1,...,x_n}$ be the law on
$\Omega_+^{x_1,...x_n}$, defined as the image of $Q$ by the application
\begin{eqnarray*}
  \Omega_+ &\to& \Omega_+ ^{x_1,...x_n} \\
  \omega &\mapsto& (\Theta_{x_i-x_1}(\omega)_+)_{1 \le i \le n}
\end{eqnarray*}
It follows from \eqref{invariance_Q} that all the one-dimensional marginals of
$Q^{x_1,...,x_n}$ are equal to $Q$. More generally, all the laws defined in that
way are compatible. Hence Kolmogorov's theorem yields the existence of a law
$\overline{Q}$ on $\Omega^\Real$ such that the finite dimensional marginal
 of $\overline{Q}$ on $x_1, ..., x_n$ is equal to $Q^{x_1,...,x_n}$, whatever $x_1< ... < x_n$.

Consider $(Z^x)_{x \in \Real}$ a variable on $\Omega^\Real$ with law $\overline Q$.
Fix $b>0$. For any $a\ge b$, we have $\Theta_{a-b}(Z^{-a})_+=Z^{-b}$, and therefore
$\Theta_{a-b}(Z^{-a})_+$ is independent of $a\ge b$.
Hence $\Theta_a(Z^{-a})$ is an element of $\Omega_+$ -- write it $\omega$ -- such that
 the value of $T_{-b}(\omega)$ and the restriction of $\omega$ to $[T_{-b},\infty)$ is independent of $a\ge b$.
We define a variable $Z$ on $\Omega_0$ by setting
 $$Z:= \lim_{a\to +\infty} \Theta_a(Z^{-a}),$$
where the limit is taken pointwise. Call $P$ its law. We clearly
have $P_+=Q$. Moreover the law $P$ is spatially stationary. Indeed,
for any $x\in \Real$, we have
$$ \Theta_x(Z) = \lim_{a \to +\infty} \Theta_{x+a} (Z^{-a}) = \lim_{a \to +\infty} \Theta_{a} (Z^{-a-x}),$$
but the family $(Z^{-a-x})_{a\in \Real}$ has the same law as
$(Z^{-a})_{a\in \Real}$, so $\Theta_{x}(Y)$ also has law $P$.
Finally, we should prove the following convergence result of
laws on $\Omega_0$:
$$P_u \circ \Theta_x \rightarrow P,$$
for any $x\in \Real$.
Take $f$ any positive bounded continuous functional depending on a
finite number of variables $\omega_{t_1}, ... \omega_{t_n}$, with
$t=t_1< ... < t_n$, so that $f((\omega_s)_{s \in
\Z})=f((\omega_s)_{s \ge t})$. We suppose without loss of generality
$t<0$. Observe that under the probability $P_u \circ \Theta_x$ or
under $P$, we have $T_0=0$, and the events $T_{-y} \le t$ and $T_y
\circ \Theta_{-y} >-t$ coincide, almost surely. Observe also
$Q(T_y\le -t) \to_{y \to \infty} 0$. Then,
 \begin{eqnarray*}
P_u\circ \Theta_x (f((\omega_{s})_{s\ge t}) \un_{T_{-y} < t}) &=&
 P_u \circ \Theta_{x-y} (f \circ \Theta_y
((\omega_{s})_{s\ge t}), T_{y} >-t) \\
&=&  (P_u \circ \Theta_{x-y})_+ (f \circ \Theta_y
((\omega_{s})_{s\ge t}), T_{y} >-t) \\
&\underset {u \to 0^+} \to& Q(f \circ \Theta_y ((\omega_{s})_{s\ge t}), T_{y} > -t) \\
&=& P(f \circ \Theta_y ((\omega_{s})_{s\ge t}), T_{y} >-t) \\
&=& P(f((\omega_{s})_{s\ge t}), T_{-y} <t),
\end{eqnarray*}
where we get the second line because the functional $\un_{T_{y} >-t} f \circ
\Theta_y ((\omega_{s})_{s\ge t}) $ does not depend on
$(\omega_n)_{n<0}$, and where we obtain the last line thanks to the translation
$\Theta_{-y}$. Besides, we have:
\begin{eqnarray*}
&& |P_u\circ \Theta_x (f((\omega_{s})_{s\ge t}) \un_{T_{-y} < t})-
P_u \circ \Theta_x (f((\omega_{s})_{s\ge t}) )| \\
&\le& (\sup f).\ P_u \circ \Theta_x (\un_{T_{-y} \ge t}) \\
&=& (\sup f).\  P_u\circ \Theta_{x-y} (\un_{T_{y} \le -t})\\
&\underset{v \to 0^+} \to&  (\sup f).\  Q(\un_{T_{y} \le -t})
\underset{y\to \infty} \to 0,
 \end{eqnarray*}
and also
$$P(f((\omega_{s})_{s\ge t}), T_{-y} <-t) \underset{y\to \infty} \To P(f((\omega_{s})_{s\ge t})).$$
This is enough to deduce
 $$P_u \circ \Theta_x (f((\omega_{s})_{s\ge t}) \underset{u\to 0^+} \To P(f((\omega_{s})_{s\ge t})).$$
The law $P_u \circ \Theta_x$ converges weakly to $P$.

\section*{Acknowledgements}
I would like to thank specially J.\ Bertoin, my supervisor during
this work, and A.\ Lachal, who made a remarkable work of rereading,
and suggested many corrections and improvements.

\bibliographystyle{abbrv}
\bibliography{biblio}
\end{document}